\documentclass[11pt]{amsart}
\usepackage{mathrsfs}
\usepackage{amssymb,latexsym}
\usepackage{pstcol,pstricks,color}

 \numberwithin{equation}{section}

\newtheorem{theorem}{Theorem}[section]

\newtheorem{corollary}[theorem]{Corollary}

\newtheorem{remark}[theorem]{Remark}

\newtheorem{lemma}[theorem]{Lemma}

\def\qed{\hfill $\Box$}
\def\pf{\noindent {\it Proof.} }

\title{Symmetric and asymmetric peaks or valleys in (partial) Dyck paths }

\begin{document}
\maketitle
\begin{center}
Yidong Sun$^{\dag}$\footnote{Corresponding author: Yidong Sun.}, Wenle Shi$^{\ddag}$ and Di Zhao$^{\sharp}$

School of Science, Dalian Maritime University, 116026 Dalian, P.R. China\\[5pt]

{\it Emails: $^{\dag}$sydmath@dlmu.edu.cn, $^{\ddag}$wlshi@dlmu.edu.cn, $^{\sharp}$zd1129@dlmu.edu.cn  }

\end{center}\vskip0.2cm

\subsection*{Abstract} The concepts of symmetric and asymmetric peaks in Dyck paths were introduced by
Fl\'{o}rez and Ram\'{\i}rez, who counted the total number of such peaks over all Dyck
paths of a given length. Elizalde generalized their results by giving multivariate
generating functions that keep track of the number of symmetric peaks and the number
of asymmetric peaks. Elizalde also considered the analogous notion of symmetric valleys by a continued fraction method. In this paper, mainly by bijective methods, we devote to
enumerating the statistics ``symmetric peaks", ``asymmetric peaks", ``symmetric valleys" and ``asymmetric valleys" of weight $k+1$ over all (partial) Dyck paths of a given length. Our results refine some consequences of
Fl\'{o}rez and Ram\'{\i}rez, and Elizalde.

\medskip

{\bf Keywords}: Dyck path, Symmetric peak, Asymmetric peak, Symmetric valley, Asymmetric valley, Riordan array.

\noindent {\sc 2010 Mathematics Subject Classification}: Primary 05A15; Secondary 05A05, 05A19.

{\bf \section{ Introduction } }

A {\it free Dyck path} of length $2n$ is a lattice path from $(0, 0)$ to $(2n, 0)$ in the $XOY$-plane and consists of up steps $\mathbf{u}=(1, 1)$ and down steps $\mathbf{d}=(1, -1)$. It is known that the set of free Dyck paths of length $2n$ is counted by $P_n=\binom{2n}{n}$, which has the generating function
\begin{eqnarray}\label{eqn 1.1}
P(x)=\sum_{n\geq 0}\binom{2n}{n}x^n=\frac{1}{\sqrt{1-4x}}.
\end{eqnarray}

A {\it Dyck path} of length $2n$ is a {\it free Dyck path} of length $2n$ that does not go below the $X$-axis. See \cite[p.204]{StanleyEC} and \cite{Deutsch99}. Let $\mathcal{D}_n$ be the set of Dyck paths of length $2n$. It is well-known \cite{Stanley} that $|\mathcal{D}_n|=C_n=\frac{1}{n+1}\binom{2n}{n}$, the $n$th Catalan number, has the generating function
\begin{eqnarray*}
C(x)=\sum_{n\geq 0}C_nx^n=\frac{1-\sqrt{1-4x}}{2x}
\end{eqnarray*}
with the relation $C(x)=1+xC(x)^2=\frac{1}{1-xC(x)}$.

A {\it partial Dyck path} of length $2n-k$ is the prefix of a Dyck path from $(0, 0)$ to $(2n-k, k)$. Let $\mathcal{D}_{n,k}$ be the set of partial Dyck paths of length $2n-k$. It is known that $|\mathcal{D}_{n,k}|=C_{n,k}=\frac{k+1}{n+1}\binom{2n-k}{n}$ has the generating function
\begin{eqnarray}\label{eqn 1.2}
x^kC(x)^{k+1}=\sum_{n\geq k}C_{n,k}x^n=\sum_{n\geq k}\frac{k+1}{n+1}\binom{2n-k}{n}x^n.
\end{eqnarray}
In fact, the martix $\big(C_{n,k}\big)_{n\geq k\geq 0}$ forms a Riordan array $(C(x), xC(x))$, the first values of $C_{n,k}$ are illustrated in Table 1.1.
\begin{center}
\begin{eqnarray*}
\begin{array}{c|cccccccc}\hline
n/k & 0   & 1   & 2    & 3    & 4    & 5    & 6     & 7   \\\hline
  0 & 1   &     &      &      &      &      &       &     \\
  1 & 1   & 1   &      &      &      &      &       &     \\
  2 & 2   & 2   & 1    &      &      &      &       &      \\
  3 & 5   & 5   & 3    & 1    &      &      &       &      \\
  4 & 14  & 14  & 9    & 4    &  1   &      &       &      \\
  5 & 42  & 42  & 28   & 14   &  5   &  1   &       &     \\
  6 & 132 & 132 & 90   & 48   &  20  &  6   &  1    &     \\
  7 & 429 & 429 & 297  & 165  &  75  &  27  &  7    & 1   \\\hline
\end{array}
\end{eqnarray*}
Table 1.1. The first values of $C_{n,k}$.
\end{center}

Recall that {\it Riordan array} \cite{ShapB, ShapGet, Sprug} is an infinite lower triangular matrix $\mathscr{D}=(d_{n,k})_{n,k \in \mathbb{N}}$ such that its $k$-th column has generating function $d(x)h(x)^k$, where $d(x)$ and $h(x)$ are formal power series with $d(0)=1$ and $h(0)=0$. That is, the general term of $\mathscr{D}$ is $d_{n,k}=[x^n]d(x)h(x)^k$, where $[x^n]$ is the coefficient operator. The matrix $\mathscr{D}$ corresponding to the pair $d(x)$ and $h(x)$ is denoted by $(d(x),h(x))$. A Riordan array $\mathscr{D}=(d(x),h(x))$ is {\it proper}, if $h'(0) \neq 0$ additionally.

Let $\varepsilon$ be the empty path, that is a dot path. If $P_1$ and $P_2$ are (partial) Dyck paths, then we define $P_1P_2$ as the concatenation of $P_1$ and $P_2$,
and define $\overline{P}_1$ as the reverse path of $P_1$. For example, $P_1=\mathbf{uuduuddd}$ and $P_2=\mathbf{uudd}$, then $P_1P_2=\mathbf{uuduuddd}\mathbf{uudd}$ and $\overline{P}_1=\mathbf{uuuddudd}$.

A point of a (partial) Dyck path with ordinate $\ell$ is said to be at {\it level} $\ell$. A step of a (partial) Dyck path is said to be at level $\ell$ if the ordinate of its endpoint is $\ell$. By a {\it return step} we mean a $\mathbf{d}$-step at level $0$. Dyck paths that have exactly one return step are said to be {\it primitive}. A {\it peak (valley)} in a (partial) Dyck path is an occurrence of $\mathbf{ud}$ ($\mathbf{du}$). By the {\it level of a peak (valley)} we mean the level of the intersection point of its two steps.
A {\it pyramid} in a (partial) Dyck path is a section of the form $\mathbf{u}^{h}\mathbf{d}^{h}$,
a succession of $h$ up steps followed immediately by $h$ down steps, where $h$ is called the {\it height} of the pyramid.
A {\it maximal mountain} of a (partial) Dyck path is a maximal subsequence of the form $\mathbf{u}^{i}\mathbf{d}^{j}$ for $i, j\geq 1$. Note that a maximal mountain contains a unique peak and vice verse.
A peak is {\it symmetric} ({\it asymmetric}) if its maximal mountain $\mathbf{u}^{i}\mathbf{d}^{j}$ satisfies $i=j$ $(i\neq j)$, and it is {\it left asymmetric} when $i>j$ and {\it right asymmetric} when $i<j$.
The {\it weight} of a peak is defined to be $\min\{i, j\}$ when its maximal mountain is $\mathbf{u}^{i}\mathbf{d}^{j}$.
The corresponding concepts to valleys of (partial) Dyck paths are defined similarly. See Figure 1 for detailed illustrations.

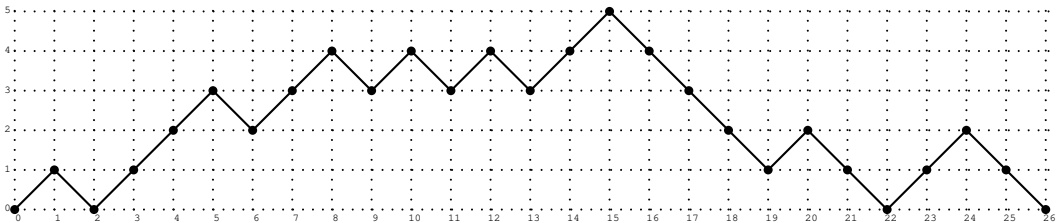
\begin{figure}[h] \setlength{\unitlength}{0.5mm}

\begin{center}
\begin{pspicture}(13.2,3)
\psset{xunit=15pt,yunit=15pt}\psgrid[subgriddiv=1,griddots=4,
gridlabels=4pt](0,0)(26,5)

\psline(0,0)(1,1)(2,0)(5,3)(6,2)(8,4)(9,3)(10,4)
\psline(10,4)(11,3)(12,4)(13,3)(15,5)(19,1)(20,2)(22,0)(24,2)(26,0)

\pscircle*(0,0){0.06}\pscircle*(1,1){0.06}\pscircle*(2,0){0.06}
\pscircle*(3,1){0.06}\pscircle*(4,2){0.06}\pscircle*(5,3){0.06}
\pscircle*(6,2){0.06}\pscircle*(7,3){0.06}\pscircle*(8,4){0.06}
\pscircle*(9,3){0.06}\pscircle*(10,4){0.06}\pscircle*(11,3){0.06}
\pscircle*(12,4){0.06}\pscircle*(13,3){0.06}\pscircle*(14,4){0.06}
\pscircle*(15,5){0.06}\pscircle*(16,4){0.06}\pscircle*(17,3){0.06}
\pscircle*(18,2){0.06}\pscircle*(19,1){0.06}\pscircle*(20,2){0.06}
\pscircle*(21,1){0.06}\pscircle*(22,0){0.06}\pscircle*(23,1){0.06}
\pscircle*(24,2){0.06}\pscircle*(25,1){0.06}\pscircle*(26,0){0.06}

\end{pspicture}
\end{center}\vskip0.5cm

\caption{\small A Dyck path of length $26$ with three symmetric peaks and two symmetric valleys of weight $1$, one symmetric peak and one symmetric valley of weight $2$, two left and one right asymmetric peaks of weight $1$, one right asymmetric peaks of weight $2$, one left and three right asymmetric valleys of weight $1$.}

\end{figure}

In the literature, there are many papers dedicated to statistics of Dyck paths (words), see \cite{Asakly, BlechBren, Czabarka1, Czabarka2, Deutsch99, Elizalde, Florez1, FloRam, Florez2, ManSap, SapTasTsi}, \cite{Sun, SunJia} and the references therein. Recently, Fl\'{o}res and Ram\'{\i}rez \cite{FloRam} find a formula for the total number, $sp(n)$, of symmetric peaks over all Dyck paths of length $2(n+1)$, as well as for the total number, $ap(n)$, of asymmetric peaks over all Dyck paths of length $2(n+3)$. Elizalde \cite{Elizalde} obtains a trivariate generating function that enumerates Dyck paths with respect to the number of symmetric peaks and the number of asymmetric peaks. His method gives a more direct derivation of the
generating function for $sp(n)$ and $ap(n)$. Namely,
\begin{eqnarray}
\sum_{n\geq 0}\mathrm{sp}(n)x^n \hskip-.22cm &=&\hskip-.22cm \frac{1}{2x}\Big(1+\frac{5x-1}{(1-x)\sqrt{1-4x}}\Big)=\frac{C(x)}{1-x}\Big(1+\frac{x}{\sqrt{1-4x}}\Big), \label{eqn 1.3}\\
\sum_{n\geq 0}\mathrm{ap}(n)x^n \hskip-.22cm &=&\hskip-.22cm \frac{1}{x^3}\Big(\frac{1-3x}{(1-x)\sqrt{1-4x}}-1\Big)= \frac{2C(x)^3}{(1-x)\sqrt{1-4x}}. \label{eqn 1.4}
\end{eqnarray}
The sequence $sp(n)$ reads $1, 3, 8, 23, 72, 240, 834, 2979, 10844, 40016, \dots$, and $ap(n)$ reads $2, 12$, $54, 222, 882, 3456, 13466, 52362,\dots$ for $n\geq 0.$

Elizalde \cite{Elizalde} also deals with the related notion of symmetric valleys, originally suggested by Deutsch \cite{Deutsch99},
phrased in terms of pairs of consecutive peaks at the same level. By a continued fraction method, he deduces a simple generating function for the total number, $sv(n)$, of symmetric valleys over all Dyck paths of length $2(n+2)$. That is
\begin{eqnarray}\label{eqn 1.5}
\sum_{n\geq 0}\mathrm{sv}(n)x^n \hskip-.22cm &=&\hskip-.22cm  \frac{2}{1-3x-4x^2+(1-x)\sqrt{1-4x}}=\frac{C(x)}{\sqrt{1-4x}}\frac{1}{1-x^2C(x)^2}.
\end{eqnarray}
The sequence $sv(n)$ reads $1, 3, 11, 40, 148, 553, 2083, \dots$, see $A014301$ in \cite{Sloane}.

In this paper, mainly by bijective methods, we enumerate the statistics ``symmetric peaks", ``asymmetric peaks", ``symmetric valleys" and ``asymmetric valleys" of weight $k+1$ over all (partial) Dyck paths of a given length.

\vskip0.5cm
\section{Symmetric and asymmetric peaks with weight $k+1$ in Dyck paths}

\subsection{Symmetric peaks with weight $k+1$ in Dyck paths}

In this subsection, we concentrate on the symmetric peaks with weight $k+1$ in Dyck paths.

Let $\mathcal{S}_{n,k}$ denote the set of Dyck paths of length $2(n+1)$ having a distinguished symmetric peak with weight $k+1$. Set $S_{n,k}=|\mathcal{S}_{n,k}|$, which is the total number of symmetric peaks with weight $k+1$ in $\mathcal{D}_{n+1}$.

\begin{lemma}\label{lemma 2.1.1}
There exists a bijection between the sets $\mathcal{S}_{n,k}$ and $\mathcal{S}_{n+j,k+j}$ for $j\geq 1$.
\end{lemma}
\pf Given a Dyck path $P\in \mathcal{S}_{n,k}$ with a distinguished symmetric peak $\mathbf{u}^{k+1}\mathbf{d}^{k+1}$, insert a pyramid $\mathbf{u}^j\mathbf{d}^j$ at the top of the distinguished symmetric peak to form the
distinguished symmetric peak $\mathbf{u}^{k+j+1}\mathbf{d}^{k+j+1}$, the resulting Dyck path $P'$ is in $\mathcal{S}_{n+j,k+j}$.

Conversely, for any $P'\in \mathcal{S}_{n+j,k+j}$ with a distinguished symmetric peak $\mathbf{u}^{k+j+1}\mathbf{d}^{k+j+1}$, remove the sub-path $\mathbf{u}^j\mathbf{d}^j$ to produce a Dyck path $P\in \mathcal{S}_{n,k}$ with the distinguished symmetric peak $\mathbf{u}^{k+1}\mathbf{d}^{k+1}$.    \qed\vskip0.2cm

Let $\mathcal{F}_{n,k}$ denote the set of pairs $(F, D)$, where $F$ are empty path or free Dyck paths starting with a $\mathbf{ud}$ and $D$ are partial Dyck paths ending at level $k$ such that the length sum of $F$ and $D$ is $2n-k$. When $k=0$, $D$ is naturally a Dyck path including the empty case. Let $F_{n, k}=|\mathcal{F}_{n,k}|$. Then by (\ref{eqn 1.1}), (\ref{eqn 1.2}) and the expansions
\begin{eqnarray}\label{eqn 2.1}
\sum_{n\geq 0}\binom{2n+r}{n}x^n &=& \frac{C(x)^r}{\sqrt{1-4x}}.
\end{eqnarray}
it is easy to deduce that
\begin{eqnarray*}
F_{k}(x)=\sum_{n\geq k}F_{n, k}x^n=x^{k}\Big(1+\frac{x}{\sqrt{1-4x}}\Big)C(x)^{k+1}
\end{eqnarray*}
and
\begin{eqnarray}\label{eqn 2.2}
F_{n, k} =[x^n]x^{k}\Big(1+\frac{x}{\sqrt{1-4x}}\Big)C(x)^{k+1}=\Big(\frac{k+1}{n+1}+\frac{n-k}{2n-k}\Big)\binom{2n-k}{n}.
\end{eqnarray}

The triangle $\mathbf{F}=\big(F_{n,k}\big)_{n\geq k\geq 0}$ forms a Riordan array $\big(C(x)(1+\frac{x}{\sqrt{1-4x}}), xC(x)\big)$. The first values of $F_{n,k}$ are exhibited in Table 2.1.
\begin{center}
\begin{eqnarray*}
\begin{array}{|c|cccccc|}\hline
n/k  & 0    & 1    & 2    & 3    & 4    & 5      \\\hline
  0  & 1    &      &      &      &      &         \\
  1  & 2    & 1    &      &      &      &         \\
  2  & 5    & 3    & 1    &      &      &         \\
  3  & 15   & 9    & 4    &  1   &      &         \\
  4  & 49   & 29   & 14   &  5   &  1   &         \\
  5  & 168  & 98   & 49   &  20  &  6   &  1      \\\hline
\end{array}
\end{eqnarray*}
Table 2.1. The first values of $F_{n,k}$.
\end{center}

\begin{theorem}\label{theom 2.1.2}
There is a bijection between the sets $\mathcal{S}_{n,0}$ and $\mathcal{F}_{n,0}$.
\end{theorem}
\pf Given a Dyck path $Q\in \mathcal{S}_{n,0}$ with a distinguished symmetric peak $\underline{\mathbf{u}\mathbf{d}}$, when $\underline{\mathbf{u}\mathbf{d}}$ is at level $1$, that is, $Q=P_1\underline{\mathbf{ud}}Q_1$, where $P_1, Q_1$ are Dyck paths. Then there are two cases to be considered. The first is that $P_1$ is an empty path, we define $\phi(Q)=(\varepsilon, Q_1)\in \mathcal{F}_{n,0}$. If not, then $P_1=\mathbf{u}P_2\mathbf{d}$, we define $\phi(Q)=(\mathbf{ud}P_2, Q_1)\in \mathcal{F}_{n,0}$. Note that $P_2$ is a free Dyck path above the line $y=-1$.

When the distinguished symmetric peak $\underline{\mathbf{u}\mathbf{d}}$ of $Q$ is at level $k\geq 2$, $Q$ can be uniquely partitioned into $Q=Q_2P_1\mathbf{d}\underline{\mathbf{ud}}\mathbf{u}P_2Q_1$, where $P_1\mathbf{d}\underline{\mathbf{ud}}\mathbf{u}P_2$ is a primitive Dyck path and $Q_1, Q_2$ are Dyck paths. Then define $\phi(Q)=(\mathbf{ud}P_2Q_2P_1, Q_1)\in \mathcal{F}_{n,0}$. Note that $P_2$ ends with a $\mathbf{d}$ step and $Q_2P_1$ begins with a $\mathbf{u}$ step, $P_2Q_2P_1$ is a free Dyck path such that the intersection of $P_2$ and $Q_2P_1$ forms a leftmost lowest valley at the level $-k$ and the intersection of $P_2Q_2$ and $P_1$ also forms a rightmost lowest valley at the level $-k$.

Conversely, the inverse map of $\phi$ is constructed as follows. For any $(P_1, Q_1)\in \mathcal{F}_{n,0}$, when $P_1$ is empty, then define $\phi^{-1}(\varepsilon, Q_1)=\underline{\mathbf{ud}}Q_1$, we get a Dyck path $Q=\underline{\mathbf{ud}}Q_1\in \mathcal{S}_{n,0}$ such that $Q$ starting with a symmetric peak $\underline{\mathbf{ud}}$. When $P_1=\mathbf{ud}P_2$ such that $P_2$ is a free Dyck path above the line $y=-1$, then define $\phi^{-1}(P_1, Q_1)=\mathbf{u}P_2\mathbf{d}\underline{\mathbf{ud}}Q_1$, we get a Dyck path $Q=\mathbf{u}P_2\mathbf{d}\underline{\mathbf{ud}}Q_1\in \mathcal{S}_{n,0}$ such that the distinguished symmetric peak $\underline{\mathbf{ud}}$ of $Q$ is  at level $1$ and not at the beginning of $Q$.
When $P_1=\mathbf{ud}P_2$ such that $P_2$ is a free Dyck path with the lowest valley at the level $-k$ for $k\geq 2$, according to the leftmost and rightmost lowest valleys (which have the same level $-k$), $P_2$ can be uniquely written as $P_2=P_3Q_2P_4$, where $Q_2$ is the Dyck path between the leftmost and rightmost lowest valleys of $P_2$. Then define $\phi^{-1}(P_1, Q_1)=Q_2P_4\mathbf{d}\underline{\mathbf{ud}}\mathbf{u}P_3Q_1$, we get a Dyck path $Q=Q_2P_4\mathbf{d}\underline{\mathbf{ud}}\mathbf{u}P_3Q_1\in \mathcal{S}_{n,0}$ such that the distinguished symmetric peak $\underline{\mathbf{ud}}$ of $Q$ is  at level $k\geq 2$.   \qed\vskip0.2cm

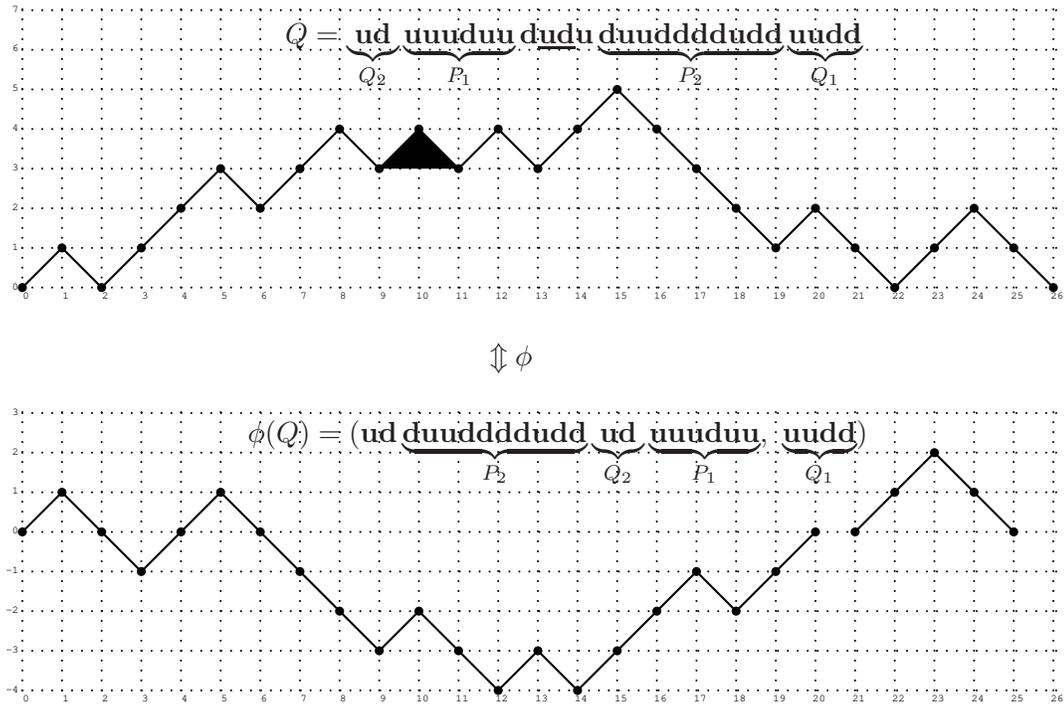
\begin{figure}[h] \setlength{\unitlength}{0.5mm}

\begin{center}
\begin{pspicture}(13,4)
\psset{xunit=15pt,yunit=15pt}\psgrid[subgriddiv=1,griddots=4,
gridlabels=4pt](0,0)(26,7)

\psline(0,0)(1,1)(2,0)(5,3)(6,2)(8,4)(9,3)(10,4)
\psline(10,4)(11,3)(12,4)(13,3)(15,5)(19,1)(20,2)(22,0)(24,2)(26,0)

\pscircle*(0,0){0.06}\pscircle*(1,1){0.06}\pscircle*(2,0){0.06}
\pscircle*(3,1){0.06}\pscircle*(4,2){0.06}\pscircle*(5,3){0.06}
\pscircle*(6,2){0.06}\pscircle*(7,3){0.06}\pscircle*(8,4){0.06}
\pscircle*(9,3){0.06}\pscircle*(10,4){0.06}\pscircle*(11,3){0.06}
\pscircle*(12,4){0.06}\pscircle*(13,3){0.06}\pscircle*(14,4){0.06}
\pscircle*(15,5){0.06}\pscircle*(16,4){0.06}\pscircle*(17,3){0.06}
\pscircle*(18,2){0.06}\pscircle*(19,1){0.06}\pscircle*(20,2){0.06}
\pscircle*(21,1){0.06}\pscircle*(22,0){0.06}\pscircle*(23,1){0.06}
\pscircle*(24,2){0.06}\pscircle*(25,1){0.06}\pscircle*(26,0){0.06}

\psline*(9,3)(10,4)(11,3)

\put(3.5,3.25){$Q=\underbrace{\mathbf{ud}}_{Q_2}\underbrace{\mathbf{uuuduu}}_{P_1}\mathbf{d}\underline{\mathbf{ud}}
\mathbf{u}\underbrace{\mathbf{duuddddudd}}_{P_2}\underbrace{\mathbf{uudd}}_{Q_1}$}

\end{pspicture}
\end{center}\vskip0.5cm

$\Updownarrow \phi$

\begin{center}
\begin{pspicture}(13,2)
\psset{xunit=15pt,yunit=15pt}\psgrid[subgriddiv=1,griddots=4,
gridlabels=4pt](0,-4)(26,3)

\psline(0,0)(1,1)(2,0)(3,-1)(5,1)(9,-3)(10,-2)(12,-4)(13,-3)(14,-4)(17,-1)(18,-2)(20,0)

\psline(21,0)(23,2)(25,0)

\pscircle*(0,0){0.06}\pscircle*(1,1){0.06}\pscircle*(2,0){0.06}
\pscircle*(3,-1){0.06}\pscircle*(4,0){0.06}\pscircle*(5,1){0.06}
\pscircle*(6,0){0.06}\pscircle*(7,-1){0.06}\pscircle*(8,-2){0.06}
\pscircle*(9,-3){0.06}\pscircle*(10,-2){0.06}\pscircle*(11,-3){0.06}
\pscircle*(12,-4){0.06}\pscircle*(13,-3){0.06}\pscircle*(14,-4){0.06}\pscircle*(15,-3){0.06}
\pscircle*(16,-2){0.06}\pscircle*(17,-1){0.06}\pscircle*(18,-2){0.06}
\pscircle*(19,-1){0.06}\pscircle*(20,0){0.06}
\pscircle*(21,0){0.06}\pscircle*(22,1){0.06}\pscircle*(23,2){0.06}
\pscircle*(24,1){0.06}\pscircle*(25,0){0.06}

\put(3,1.2){$\phi(Q)=(\mathbf{ud}\underbrace{\mathbf{duuddddudd}}_{P_2}\underbrace{\mathbf{ud}}_{Q_2}\underbrace{\mathbf{uuuduu}}_{P_1},\  \underbrace{\mathbf{uudd}}_{Q_1})$}

\end{pspicture}
\end{center}
\vskip2cm

\caption{\small An example of the bijection $\phi$ described in the proof of Theorem \ref{theom 2.1.2}. }

\end{figure}

In order to give a more intuitive view on the bijection $\phi$, we present a pictorial description of
$\phi$ for the case $Q=\mathbf{uduuuduu}\mathbf{d}\underline{\mathbf{ud}}\mathbf{u}\mathbf{duuddddudduudd}\in \mathcal{S}_{12,0}$ and $\phi(Q)=(\mathbf{ud}\mathbf{duuddddudd}\mathbf{ud}\mathbf{uuuduu}, \mathbf{uudd})\in \mathcal{F}_{12,0}$. See Figure 2 for detailed illustrations.

Let $S_k(x)=\sum_{n\geq k}S_{n,k}x^n$, then $S_0(x)=F_0(x)$ by Theorem \ref{theom 2.1.2}, and $S_k(x)=x^kS_0(x)$ by Lemma \ref{lemma 2.1.1}. Together with (\ref{eqn 2.2}), we have

\begin{corollary}\label{coro 2.1.3} The generating function
\begin{eqnarray*}
S_k(x)=\sum_{n\geq k}S_{n,k}x^n = x^kC(x)\Big(1+\frac{x}{\sqrt{1-4x}}\Big),
\end{eqnarray*}
and the triangle $\mathbf{S}=\big(S_{n,k}\big)_{n\geq k\geq 0}$ forms a Riordan array $\big(C(x)(1+\frac{x}{\sqrt{1-4x}}), x\big)$ with the general entry
$$S_{n,k}=\frac{n-k+3}{2}C_{n-k}$$
for $n>k$ and $S_{n,n}=1$ for $n\geq 0$.
\end{corollary}

The first values of $S_{n,k}$ are showed in Table 2.2.

\begin{center}
\begin{eqnarray*}
\begin{array}{|c|cccccc|}\hline
n/k & 0   & 1    & 2    & 3    & 4    & 5       \\\hline
  0 & 1   &      &      &      &      &        \\
  1 & 2   & 1    &      &      &      &         \\
  2 & 5   & 2    & 1    &      &      &         \\
  3 & 15  & 5    & 2    &  1   &      &         \\
  4 & 49  & 15   & 5    &  2   &  1   &        \\
  5 & 168 & 49   & 15   &  5   &  2   &  1     \\\hline
\end{array}
\end{eqnarray*}
Table 2.2. The first values of $S_{n,k}$.
\end{center}

Naturally, the total number $sp(n)$ of symmetric peaks is the row sum of the triangle $\mathbf{S}$, i.e.,
\begin{eqnarray*}
\mathrm{sp}(n)=1+\sum_{i=1}^{n}\frac{i+3}{2}C_{i}
\end{eqnarray*}
and has the generating function given by (\ref{eqn 1.3}).

\subsection{Asymmetric peaks with weight $k+1$ in Dyck paths}

In this subsection, we consider the left asymmetric peaks with weight $k+1$ in Dyck paths. The right asymmetric peaks are equivalent distribution to the left asymmetric peaks according to the symmetry of Dyck paths.

Let $\mathcal{L}_{n,k}$ denote the set of Dyck paths of length $2(n+3)$ having a distinguished left asymmetric peak with weight $k+1$. Set $L_{n,k}=|\mathcal{L}_{n,k}|$, which is the total number of left asymmetric peaks with weight $k+1$ in $\mathcal{D}_{n+3}$.

\begin{lemma}\label{lemma 2.2.1}
There exists a bijection between the sets $\mathcal{L}_{n,k}$ and $\mathcal{L}_{n+j,k+j}$ for $j\geq 1$.
\end{lemma}
\pf Similar to the proof of Lemma \ref{lemma 2.1.1}, the bijection can be constructed if one notices that a distinguished left asymmetric peak $\mathbf{u}^{k+i+1}\mathbf{d}^{k+1}$ of $P\in \mathcal{L}_{n,k}$ for certain $i\geq 1$ can be extended to a distinguished left asymmetric peak $\mathbf{u}^{k+i+j+1}\mathbf{d}^{k+j+1}$ of $Q\in \mathcal{L}_{n+j,k+j}$ by inserting a pyramid $\mathbf{u}^j\mathbf{d}^j$ at the top of $\mathbf{u}^{k+i+1}\mathbf{d}^{k+1}$, and vice verse. \qed\vskip0.2cm

Let $\mathcal{E}_{n,k}$ denote the set of pairs $(F, D)$, where $F$ are free Dyck paths and $D$ are partial Dyck paths ending at the level $k$ such that the length sum of $F$ and $D$ is $2n-k$.
Let $E_{n, k}=|\mathcal{E}_{n,k}|$. Then by (\ref{eqn 1.1}), (\ref{eqn 1.2}) and (\ref{eqn 2.2})
it is easy to deduce that
\begin{eqnarray*}
E_{k}(x)=\sum_{n\geq k}E_{n, k}x^n=\frac{x^{k}C(x)^{k+1}}{\sqrt{1-4x}}
\end{eqnarray*}
and
\begin{eqnarray}\label{eqn 2.3}
E_{n,k} =[x^n]\frac{x^{k}C(x)^{k+1}}{\sqrt{1-4x}}=\binom{2n-k+1}{n-k}.
\end{eqnarray}

The triangle $\mathbf{E}=\big(E_{n,k}\big)_{n\geq k\geq 0}$ forms a Riordan array $\big(\frac{C(x)}{\sqrt{1-4x}}, xC(x)\big)$. The first values of $E_{n,k}$ are displayed in Table 2.3.
\begin{center}
\begin{eqnarray*}
\begin{array}{|c|cccccc|}\hline
n/k  & 0    & 1    & 2    & 3    & 4    & 5      \\\hline
  0  & 1    &      &      &      &      &         \\
  1  & 3    & 1    &      &      &      &         \\
  2  & 10   & 4    & 1    &      &      &         \\
  3  & 35   & 15   & 5    &  1   &      &         \\
  4  & 126  & 56   & 21   &  6   &  1   &         \\
  5  & 462  & 210  & 84   &  28  &  7   &  1      \\\hline
\end{array}
\end{eqnarray*}
Table 2.3. The first values of $E_{n,k}$.
\end{center}

\begin{theorem}\label{theom 2.2.2}
There is a bijection between the sets $\mathcal{L}_{n,0}$ and $\mathcal{E}_{n+2,2}$.
\end{theorem}

\pf Given a Dyck path $Q\in \mathcal{L}_{n,0}$ with a distinguished left asymmetric peak $\underline{\mathbf{u}^{j+2}\mathbf{d}}$ at level $i+j+2$ for certain $i, j\geq 0$, $Q$ can be uniquely partitioned into
\begin{eqnarray*}
Q=\left\{
\begin{array}{rl}
Q_1\underline{\mathbf{u}^{j+2}\mathbf{d}}\mathbf{u}Q_3\mathbf{d}Q_4\mathbf{d}Q_5, & \mbox{when}\  i=j=0, \\[5pt]
Q_1\underline{\mathbf{u}^{j+2}\mathbf{d}}\mathbf{u}Q_3\mathbf{d}Q_4\mathbf{d}Q_5\mathbf{d}Q_2, & \mbox{when}\  i+j\geq 1,
\end{array}\right.
\end{eqnarray*}
where $Q_3, Q_4$ and $Q_5$ are Dyck paths, $Q_1$ is empty or a nonempty partial Dyck path ending with a $\mathbf{d}$ step at level $i$, and $\overline{Q}_2$ is a partial Dyck path ending at level $i+j-1\geq 0$.

In the $i=j=0$ case, $Q_1$ is always a Dyck path, we define $\varphi(Q)=(Q_1, Q_5\mathbf{u}Q_4\mathbf{u}Q_3)\in \mathcal{E}_{n+2,2}$.

In the $i+j\geq 1$ case, we define $\varphi(Q)=(Q_2\mathbf{d}Q_1\mathbf{u}^j, Q_5\mathbf{u}Q_4\mathbf{u}Q_3)\in \mathcal{E}_{n+2,2}$. Note that $Q_1$ always begins with a $\mathbf{u}$ step and ends with a $\mathbf{d}$ step if it is not empty, and $Q_2\mathbf{d}Q_1\mathbf{u}^j$ is a free Dyck path with lowest valleys at the level $-(i+j)\leq -1$ such that the leftmost lowest valley is the intersection of $Q_2\mathbf{d}$ and $Q_1\mathbf{u}^j$.

Conversely, the inverse map of $\varphi$ is constructed as follows. For any $(P_1, P_2)\in \mathcal{E}_{n+2,2}$, where $P_2=P_5\mathbf{u}P_4\mathbf{u}P_3$ and $P_3, P_4, P_5$ are Dyck paths.
When $P_1$ is a Dyck path, then define $\varphi^{-1}(P_1, P_2)=P_1\underline{\mathbf{u^2d}}\mathbf{u}P_3\mathbf{d}P_4\mathbf{d}P_5$, we get a Dyck path $Q=P_1\underline{\mathbf{u^2d}}\mathbf{u}P_3\mathbf{d}P_4\mathbf{d}P_5\in \mathcal{L}_{n,0}$ such that $Q$ has a distinguished left asymmetric peak $\underline{\mathbf{u}^{2}\mathbf{d}}$ at level $2$.

If $P_1$ is a nonempty free Dyck path with the lowest valley at the level $-k$ for $k\geq 1$, according to the leftmost lowest valley and the last maximal subpath $\mathbf{u}^j$ of $P_1$, $P_1$ can be uniquely written as $P_1=Q_2\mathbf{d}Q_1\mathbf{u}^{j}$ for certain $j\geq 0$, where the intersection of $Q_2\mathbf{d}$ and $Q_1\mathbf{u}^j$ forms the leftmost lowest valley of $P_1$.
Note that $Q_1$ is a partial Dyck path ending at level $i$ for certain $i\geq 0$ and $\overline{Q}_2$ is a partial Dyck path ending at level $i+j-1$, once $i=0$, then $j\geq 1$. The maximality of the subpath $\mathbf{u}^j$ implies that once $Q_1$ is not empty, then it ends with a $\mathbf{d}$ step. Then define $\varphi^{-1}(P_1, P_2)=Q=Q_1\underline{\mathbf{u}^{j+2}\mathbf{d}}\mathbf{u}P_3\mathbf{d}P_4\mathbf{d}P_5\mathbf{d}Q_2$, we get a Dyck path $Q\in \mathcal{L}_{n,0}$ such that the distinguished symmetric peak $\underline{\mathbf{u}^{j+2}\mathbf{d}}$ of $Q$ is at level $i+j+2\geq 3$.   \qed\vskip0.2cm

In order to give a more intuitive view on the bijection $\varphi$, we present a pictorial description of
$\varphi$ for the case $Q=\mathbf{uduuud}\underline{\mathbf{u^2d}}\mathbf{u}\mathbf{ud}\mathbf{d}\mathbf{udud}\mathbf{ddudduudd}\in \mathcal{L}_{10,0}$ and $\varphi(Q)=(\mathbf{ud}\mathbf{duuddduduuud}, \mathbf{uududuud})\in \mathcal{E}_{12,2}$. See Figure 3 for detailed illustrations.

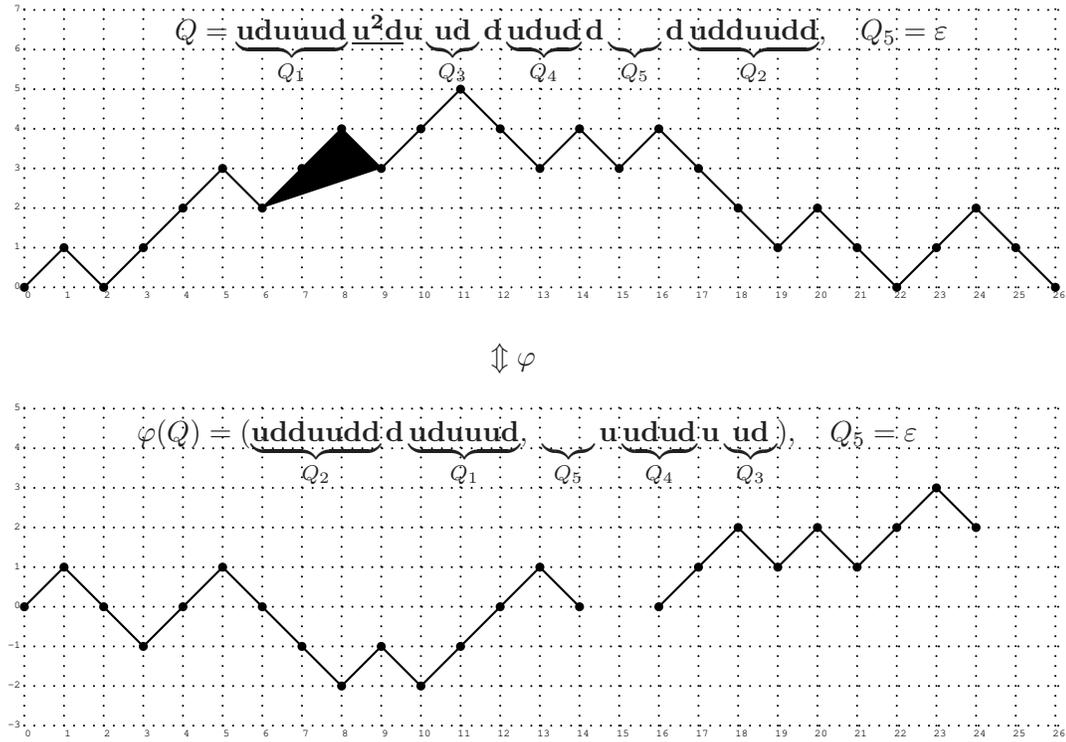
\begin{figure}[h] \setlength{\unitlength}{0.5mm}
\begin{center}
\begin{pspicture}(13,4)
\psset{xunit=15pt,yunit=15pt}\psgrid[subgriddiv=1,griddots=4,
gridlabels=4pt](0,0)(26,7)

\psline(0,0)(1,1)(2,0)(5,3)(6,2)(8,4)(9,3)(10,4)
\psline(10,4)(11,5)(12,4)(13,3)(14,4)(15,3)(16,4)(19,1)(20,2)(22,0)(24,2)(26,0)

\pscircle*(0,0){0.06}\pscircle*(1,1){0.06}\pscircle*(2,0){0.06}
\pscircle*(3,1){0.06}\pscircle*(4,2){0.06}\pscircle*(5,3){0.06}
\pscircle*(6,2){0.06}\pscircle*(7,3){0.06}\pscircle*(8,4){0.06}
\pscircle*(9,3){0.06}\pscircle*(10,4){0.06}\pscircle*(11,5){0.06}
\pscircle*(12,4){0.06}\pscircle*(13,3){0.06}\pscircle*(14,4){0.06}
\pscircle*(15,3){0.06}\pscircle*(16,4){0.06}\pscircle*(17,3){0.06}
\pscircle*(18,2){0.06}\pscircle*(19,1){0.06}\pscircle*(20,2){0.06}
\pscircle*(21,1){0.06}\pscircle*(22,0){0.06}\pscircle*(23,1){0.06}
\pscircle*(24,2){0.06}\pscircle*(25,1){0.06}\pscircle*(26,0){0.06}

\put(2,3.3){$Q=\underbrace{\mathbf{uduuud}}_{Q_1}\underline{\mathbf{u^2d}}\mathbf{u}\underbrace{\mathbf{ud}}_{Q_3}\mathbf{d}
\underbrace{\mathbf{udud}}_{Q_4}\mathbf{d}\underbrace{\empty}_{Q_5}\mathbf{d}\underbrace{\mathbf{udduudd}}_{Q_2}, \ \ \ Q_5=\varepsilon$}

\psline*(6,2)(8,4)(9,3)

\end{pspicture}
\end{center}\vskip0.5cm

$\Updownarrow \varphi$

\begin{center}
\begin{pspicture}(13,3)
\psset{xunit=15pt,yunit=15pt}\psgrid[subgriddiv=1,griddots=4,
gridlabels=4pt](0,-3)(26,5)

\psline(0,0)(1,1)(2,0)(3,-1)(5,1)(8,-2)(9,-1)(10,-2)(13,1)(14,0)

\psline(16,0)(18,2)(19,1)(20,2)(21,1)(23,3)(24,2)

\pscircle*(0,0){0.06}\pscircle*(1,1){0.06}\pscircle*(2,0){0.06}
\pscircle*(3,-1){0.06}\pscircle*(4,0){0.06}\pscircle*(5,1){0.06}
\pscircle*(6,0){0.06}\pscircle*(7,-1){0.06}\pscircle*(8,-2){0.06}
\pscircle*(9,-1){0.06}\pscircle*(10,-2){0.06}\pscircle*(11,-1){0.06}
\pscircle*(12,0){0.06}\pscircle*(13,1){0.06}\pscircle*(14,0){0.06}
\pscircle*(16,0){0.06}\pscircle*(17,1){0.06}\pscircle*(18,2){0.06}
\pscircle*(19,1){0.06}\pscircle*(20,2){0.06}
\pscircle*(21,1){0.06}\pscircle*(22,2){0.06}\pscircle*(23,3){0.06}
\pscircle*(24,2){0.06}

\put(1.5,2.2){$\varphi(Q)=(\underbrace{\mathbf{udduudd}}_{Q_2}\mathbf{d}\underbrace{\mathbf{uduuud}}_{Q_1},\  \underbrace{\empty}_{Q_5}\mathbf{u}\underbrace{\mathbf{udud}}_{Q_4}\mathbf{u}\underbrace{\mathbf{ud}}_{Q_3}), \ \ \   Q_5=\varepsilon$}

\end{pspicture}
\end{center}
\vskip1.5cm

\caption{\small An example of the bijection $\varphi$ described in the proof of Theorem \ref{theom 2.2.2}. }

\end{figure}

Let $L_k(x)=\sum_{n\geq k}L_{n,k}x^n$, then $L_0(x)=E_0(x)$ by Theorem \ref{theom 2.2.2}, and $L_k(x)=x^kL_0(x)$ by Lemma \ref{lemma 2.2.1}. Together with (\ref{eqn 2.3}), we have

\begin{corollary}\label{coro 2.2.3} The generating function
\begin{eqnarray*}
L_k(x)=\sum_{n\geq k}L_{n,k}x^n = \frac{x^kC(x)^3}{\sqrt{1-4x}},
\end{eqnarray*}
and the triangle $\mathbf{L}=\big(L_{n,k}\big)_{n\geq k\geq 0}$ forms a Riordan array $\big(\frac{C(x)^3}{\sqrt{1-4x}}, x\big)$ with the general entry
$$L_{n,k}=\binom{2n-2k+3}{n-k}.$$
\end{corollary}

The first values of $L_{n,k}$ are demonstrated in Table 2.4.

\begin{center}
\begin{eqnarray*}
\begin{array}{|c|cccccc|}\hline
n/k & 0    & 1    & 2    & 3    & 4    & 5           \\\hline
  0 & 1    &      &      &      &      &             \\
  1 & 5    & 1    &      &      &      &             \\
  2 & 21   & 5    & 1    &      &      &             \\
  3 & 84   & 21   & 5    & 1    &      &             \\
  4 & 330  & 84   & 21   & 5    &  1   &             \\
  5 & 1287 & 330  & 84   & 21   &  5  &  1           \\\hline
\end{array}
\end{eqnarray*}
Table 2.4. The first values of $L_{n,k}$.
\end{center}

Clearly, the total number $ap(n)$ of asymmetric peaks is twice the row sum of $\mathbf{L}$, i.e.,
\begin{eqnarray*}
\mathrm{ap}(n)=2\sum_{i=0}^{n}\binom{2n-2i+3}{n-i}=2\sum_{i=0}^{n}\binom{2i+3}{i}
\end{eqnarray*}
and has the generating function given by (\ref{eqn 1.4}).

\begin{remark}
The sequence $L_{n,0}=\binom{2n+3}{n}$ also counts many statistics in Dyck paths, see comments in \cite[A002054]{Sloane}. For examples, the total number of valleys,
number of $\mathbf{uu}'s$, number of peaks at level higher than one in all Dyck paths of length $2(n+2)$, and number of nonempty Dyck subpaths in all Dyck paths of length $2(n+1)$.
\end{remark}

Manes, Sapounakis, Tasoulas and Tsikouras \cite{ManSap} study {\it nonleft peaks} in Dyck paths, which are peaks such that the $\mathbf{u}$ steps preceding them are greater than or equal
to the $\mathbf{d}$ steps following them. They present a combinatorial construction of the set of Dyck paths of fixed length and number of nonleft peaks, and obtain
various results on the enumeration of several kinds of peaks. In our notations, a nonleft peak is just a symmetric peak or a left asymmetric peak.
Let $\mathcal{S}^{*}_{n,k}$ denote the set of Dyck paths of length $2(n+1)$ having a distinguished symmetric or left asymmetric peak with weight $k+1$. Set ${S}^{*}_{n,k}=|\mathcal{S}^{*}_{n,k}|$. Note that
$\mathcal{S}^{*}_{n,k}=\mathcal{S}_{n,k}\cup \mathcal{L}_{n-2,k}$ and ${S}^{*}_{n,k}=S_{n,k}+L_{n-2,k}$. By Lemma \ref{lemma 2.1.1} and \ref{lemma 2.2.1}, together with Theorem \ref{theom 2.1.2} and \ref{theom 2.2.2}, we have
\begin{corollary}\label{coro 2.2.4}
There exists a bijection between the sets $\mathcal{S}^{*}_{n,k}$ and $\mathcal{S}^{*}_{n+j,k+j}$ for $j\geq 1$. And the generating function
\begin{eqnarray*}
S^{*}_k(x)=\sum_{n\geq k}S^{*}_{n,k}x^n =S_k(x)+x^2L_k(x)=\frac{x^k}{\sqrt{1-4x}},
\end{eqnarray*}
and the triangle $\mathbf{S}^{*}=\big(S^{*}_{n,k}\big)_{n\geq k\geq 0}$ forms a Riordan array $\big(\frac{1}{\sqrt{1-4x}}, x\big)$ with the general entry
$${S}^{*}_{n,k}=\binom{2n-2k}{n-k}.$$
\end{corollary}

The first values of $S^{*}_{n,k}$ are presented in Table 2.5.

\begin{center}
\begin{eqnarray*}
\begin{array}{|c|cccccc|}\hline
n/k & 0    & 1    & 2    & 3    & 4    & 5           \\\hline
  0 & 1    &      &      &      &      &             \\
  1 & 2    & 1    &      &      &      &             \\
  2 & 6    & 2    & 1    &      &      &             \\
  3 & 20   & 6    & 2    & 1    &      &             \\
  4 & 70   & 20   & 6    & 2    &  1   &             \\
  5 & 252  & 70   & 20   & 6    &  2  &  1           \\\hline
\end{array}
\end{eqnarray*}
Table 2.5. The first values of $S^{*}_{n,k}$.
\end{center}

\vskip0.5cm

\section{Symmetric and asymmetric valleys with weight $k+1$ in Dyck paths}

\subsection{Symmetric valleys with weight $k+1$ in Dyck paths}
In this subsection, we study the symmetric valleys with weight $k+1$ in Dyck paths.

Let $\mathcal{V}_{n,k}$ denote the set of Dyck paths of length $2(n+2)$ having a distinguished symmetric valley with weight $k+1$. Set $V_{n,k}=|\mathcal{V}_{n,k}|$, which is the total number of symmetric valleys with weight $k+1$ in $\mathcal{D}_{n+2}$.

\begin{theorem}\label{theom 3.1.1}
There is a bijection between the sets $\mathcal{V}_{n,k}$ and $\mathcal{E}_{n,2k}$.
\end{theorem}

\pf Given a Dyck path $Q\in \mathcal{V}_{n,k}$ with a distinguished symmetric valley $\underline{\mathbf{d}^{k+1}\mathbf{u}^{k+1}}$ at level $i\geq 0$, $Q$ can be uniquely partitioned into
\begin{eqnarray*}
Q=\left\{
\begin{array}{rl}
Q_0\mathbf{u}Q_1\dots \mathbf{u}Q_{k}\mathbf{u}\underline{\mathbf{d}^{k+1}\mathbf{u}^{k+1}}\mathbf{d}Q_{k+1}\dots\mathbf{d}Q_{2k+1}, & \mbox{when}\  i=0, \\[5pt]
Q_0\mathbf{u}Q_1\dots \mathbf{u}Q_{k}\mathbf{u}\underline{\mathbf{d}^{k+1}\mathbf{u}^{k+1}}\mathbf{d}Q_{k+1}\dots\mathbf{d}Q_{2k+1}\mathbf{d}{Q'}_0, & \mbox{when}\  i\geq 1,
\end{array}\right.
\end{eqnarray*}
where $Q_1, \dots, Q_{2k+1}$ are Dyck paths, $Q_0$ is a partial Dyck path ending at level $i$, and $\overline{Q'}_0$ is a partial Dyck path ending at level $i-1\geq 0$.

In the $i=0$ case, $Q_0$ is always a Dyck path, we define $\theta(Q)=(Q_0, Q_1\mathbf{u}Q_{2}\dots \mathbf{u}Q_{2k+1})\in \mathcal{E}_{n,2k}$.

In the $i\geq 1$ case, we define $\theta(Q)=(Q_0'\mathbf{d}Q_0, Q_1\mathbf{u}Q_{2}\dots \mathbf{u}Q_{2k+1})\in \mathcal{E}_{n,2k}$. Note that $Q_0$ always begins with a $\mathbf{u}$ step, and $Q_0'\mathbf{d}Q_0$ is a free Dyck path with a lowest valleys at the level $-i\leq -1$ such that the leftmost lowest valley is the intersection of $Q_0'\mathbf{d}$ and $Q_0$.

Conversely, the inverse map of $\theta$ is constructed as follows. For any $(P_0, P)\in \mathcal{E}_{n,2k}$, where $P=P_1\mathbf{u}P_{2}\dots \mathbf{u}P_{2k+1}$ and $P_1, \dots, P_{2k+1}$ are Dyck paths.
When $P_0$ is a Dyck path, then define $\theta^{-1}(P_0, P)=Q=P_0\mathbf{u}P_1\dots \mathbf{u}P_{k}\mathbf{u}\underline{\mathbf{d}^{k+1}\mathbf{u}^{k+1}}\mathbf{d}P_{k+1}\dots\mathbf{d}P_{2k+1}$, we get a Dyck path $Q\in \mathcal{V}_{n,k}$ such that $Q$ has a distinguished symmetric valley $\underline{\mathbf{d}^{k+1}\mathbf{u}^{k+1}}$ at level $0$.

If $P_0$ is a nonempty free Dyck path with the lowest valley at the level $-i$ for $i\geq 1$, according to the leftmost lowest valley of $P_0$, $P_0$ can be uniquely written as $P_0={Q'}_0\mathbf{d}Q_0$, where the intersection of ${Q'}_0\mathbf{d}$ and $Q_0$ forms the leftmost lowest valley. Note that $Q_0$ and $\overline{Q'}_0$ are partial Dyck paths ending at level $i$ and $i-1$ respectively.
Then define $\theta^{-1}(P_0, P)=Q=Q_0\mathbf{u}P_1\dots \mathbf{u}P_{k}\mathbf{u}\underline{\mathbf{d}^{k+1}\mathbf{u}^{k+1}}\mathbf{d}P_{k+1}\dots\mathbf{d}P_{2k+1}\mathbf{d}{Q'}_0$, we get a Dyck path $Q\in \mathcal{V}_{n,k}$ such that the distinguished symmetric valley $\underline{\mathbf{d}^{k+1}\mathbf{u}^{k+1}}$ of $Q$ is at level $i\geq 1$.   \qed\vskip0.2cm

In order to give a more intuitive view on the bijection $\theta$, we present a pictorial description of
$\theta$ for the case $Q=\mathbf{udu^3du^3du}\underline{\mathbf{d^2u^2}}\mathbf{d^2ud^4u^2d^2}\in \mathcal{V}_{11,1}$ and $\theta(Q)=(\mathbf{d^2u^2d^3udu^3du}, \mathbf{udu^3d})\in \mathcal{E}_{11,2}$. See Figure 4 for detailed illustrations.

\begin{figure}[h] \setlength{\unitlength}{0.5mm}

\begin{center}
\begin{pspicture}(13,4)
\psset{xunit=15pt,yunit=15pt}\psgrid[subgriddiv=1,griddots=4,
gridlabels=4pt](0,0)(26,7)

\psline(0,0)(1,1)(2,0)(5,3)(6,2)(8,4)(9,5)(10,4)
\psline(10,4)(11,5)(12,4)(13,3)(14,4)(15,5)(16,4)(17,3)(18,4)(20,2)(22,0)(24,2)(26,0)

\pscircle*(0,0){0.06}\pscircle*(1,1){0.06}\pscircle*(2,0){0.06}
\pscircle*(3,1){0.06}\pscircle*(4,2){0.06}\pscircle*(5,3){0.06}
\pscircle*(6,2){0.06}\pscircle*(7,3){0.06}\pscircle*(8,4){0.06}
\pscircle*(9,5){0.06}\pscircle*(10,4){0.06}\pscircle*(11,5){0.06}
\pscircle*(12,4){0.06}\pscircle*(13,3){0.06}\pscircle*(14,4){0.06}
\pscircle*(15,5){0.06}\pscircle*(16,4){0.06}\pscircle*(17,3){0.06}
\pscircle*(18,4){0.06}\pscircle*(19,3){0.06}\pscircle*(20,2){0.06}
\pscircle*(21,1){0.06}\pscircle*(22,0){0.06}\pscircle*(23,1){0.06}
\pscircle*(24,2){0.06}\pscircle*(25,1){0.06}\pscircle*(26,0){0.06}

\put(2,3.3){$Q=\underbrace{\mathbf{uduuudu}}_{Q_0}\mathbf{u}\underbrace{\mathbf{ud}}_{Q_1}\mathbf{u}\underline{\mathbf{d^2u^2}}\mathbf{d}
\underbrace{\empty}_{Q_2}\mathbf{d}\underbrace{\mathbf{ud}}_{Q_3}\mathbf{d}\underbrace{\mathbf{dduudd}}_{Q_0'}, \ \ \ Q_2=\varepsilon$}

\psline*(11,5)(13,3)(15,5)

\end{pspicture}
\end{center}\vskip0.5cm

$\Updownarrow \theta$

\begin{center}
\begin{pspicture}(12,3)
\psset{xunit=15pt,yunit=15pt}\psgrid[subgriddiv=1,griddots=4,
gridlabels=4pt](0,-3)(24,5)

\psline(0,0)(2,-2)(4,0)(6,-2)(7,-3)(8,-2)(9,-3)(10,-2)(12,0)(13,-1)(14,0)

\psline(16,0)(17,1)(18,0)(19,1)(20,2)(21,3)(22,2)

\pscircle*(0,0){0.06}\pscircle*(1,-1){0.06}\pscircle*(2,-2){0.06}
\pscircle*(3,-1){0.06}\pscircle*(4,0){0.06}\pscircle*(5,-1){0.06}
\pscircle*(6,-2){0.06}\pscircle*(7,-3){0.06}\pscircle*(8,-2){0.06}
\pscircle*(9,-3){0.06}\pscircle*(10,-2){0.06}\pscircle*(11,-1){0.06}
\pscircle*(12,0){0.06}\pscircle*(13,-1){0.06}\pscircle*(14,0){0.06}
\pscircle*(16,0){0.06}\pscircle*(17,1){0.06}\pscircle*(18,0){0.06}
\pscircle*(19,1){0.06}\pscircle*(20,2){0.06}
\pscircle*(21,3){0.06}\pscircle*(22,2){0.06}

\put(1.5,2.2){$\theta(Q)=(\underbrace{\mathbf{dduudd}}_{Q_0'}\mathbf{d}\underbrace{\mathbf{uduuudu}}_{Q_0},\  \underbrace{\mathbf{ud}}_{Q_1}\mathbf{u}\underbrace{\empty}_{Q_2}\mathbf{u}\underbrace{\mathbf{ud}}_{Q_3}), \ \ \   Q_2=\varepsilon$}

\end{pspicture}
\end{center}
\vskip1.5cm

\caption{\small An example of the bijection $\theta$ described in the proof of Theorem \ref{theom 3.1.1}. }

\end{figure}
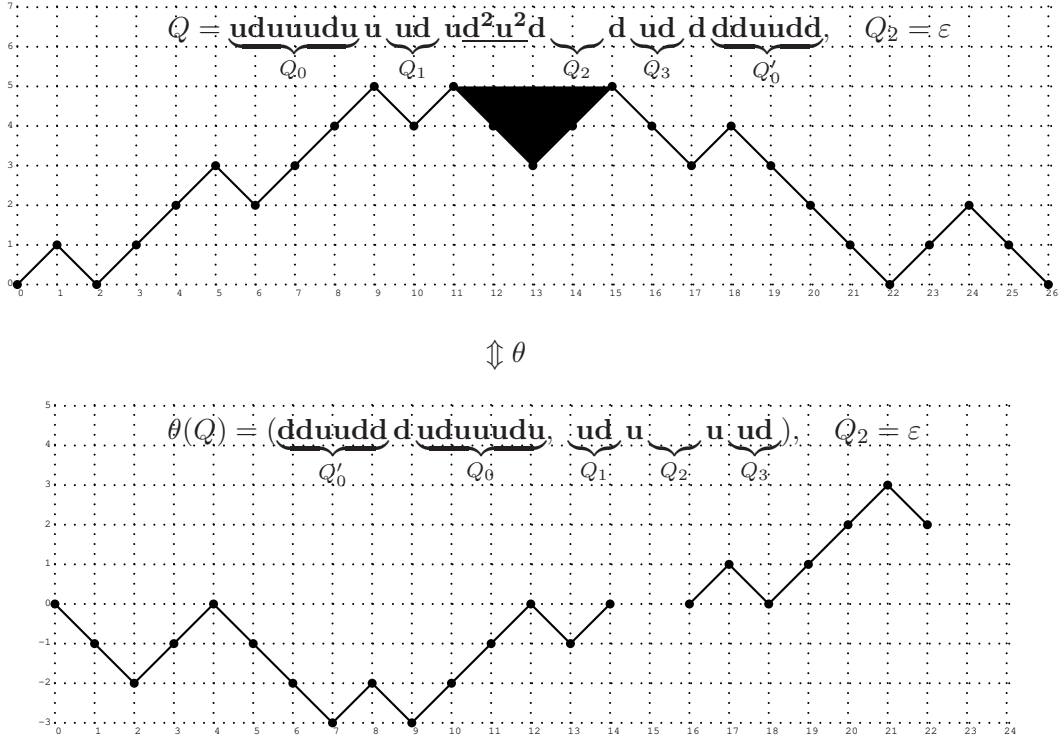

Let $V_k(x)=\sum_{n\geq k}V_{n,k}x^n$, then $V_k(x)=E_{2k}(x)$ by Theorem \ref{theom 3.1.1}. Together with (\ref{eqn 2.3}), we have

\begin{corollary}\label{coro 3.1.2} The generating function
\begin{eqnarray*}
V_k(x)=\sum_{n\geq k}V_{n,k}x^n = \frac{x^{2k}C(x)^{2k+1}}{\sqrt{1-4x}},
\end{eqnarray*}
and the triangle $\mathbf{V}=\big(V_{n,k}\big)_{n\geq k\geq 0}$ forms a Riordan array $\big(\frac{C(x)}{\sqrt{1-4x}}, x^2C(x)^2\big)$ with the general entry
$$V_{n,k}=\binom{2n-2k+1}{n-2k}.$$
\end{corollary}
The first values of $V_{n,k}$ are exhibited in Table 3.1.

\begin{center}
\begin{eqnarray*}
\begin{array}{|c|cccc|}\hline
n/k & 0    & 1      & 2    & 3              \\\hline
  0 & 1    &        &      &               \\
  1 & 3    &        &      &               \\
  2 & 10   & 1      &      &               \\
  3 & 35   & 5      &      &               \\
  4 & 126  & 21     & 1    &               \\
  5 & 462  & 84     & 7    &                \\
  6 & 1716 & 330    & 36   &  1             \\ \hline
\end{array}
\end{eqnarray*}
Table 3.1. The first values of $V_{n,k}$.
\end{center}

Obviously, the total number $sv(n)$ of symmetric valleys is the row sum of the triangle $\mathbf{V}$, i.e.,
\begin{eqnarray*}
\mathrm{sv}(n)=\sum_{k=0}^{n}\binom{2n-2k+1}{n-2k}=\sum_{k=0}^{n}\binom{2k+1}{n+1}
\end{eqnarray*}
and has the generating function given by (\ref{eqn 1.5}).

Theorem \ref{theom 2.2.2} and Theorem \ref{theom 3.1.1} in the $k=1$ case suggest the following result.

\begin{corollary}\label{coro 3.1.3}
There is a bijection between the sets $\mathcal{L}_{n,0}$ and $\mathcal{V}_{n+2,1}$.
\end{corollary}
Here we provide a simple and direct bijection. Given a Dyck path $Q\in \mathcal{L}_{n,0}$ with a distinguished left asymmetric peak $\underline{\mathbf{u}^{j+2}\mathbf{d}}$ at level $i+j+2$ for certain $i, j\geq 0$, $Q$ can be uniquely partitioned into $Q=Q_1\underline{\mathbf{u}^{j+2}\mathbf{d}}\mathbf{u}Q_2\mathbf{d}Q_3$, where $Q_1$ is empty, or a nonempty partial Dyck path ending with a $\mathbf{d}$ step at level $i$, $Q_2$ is a Dyck path and $\overline{Q}_3$ is a partial Dyck path ending at level $i+j+1\geq 1$. Then define $\eta(Q)=Q_1\mathbf{u}^{j+1}Q_2\mathbf{u}\underline{\mathbf{d}^2\mathbf{u}^2}\mathbf{d}Q_3$, we get a Dyck path $P=Q_1\mathbf{u}^{j+1}Q_2\mathbf{u}\underline{\mathbf{d}^2\mathbf{u}^2}\mathbf{d}Q_3 \in \mathcal{V}_{n+2,1}$ with a distinguished symmetric valley $\underline{\mathbf{d}^2\mathbf{u}^2}$ at level $i+j$. It is easily to verify that $\eta$ is a bijection between the sets $\mathcal{L}_{n,0}$ and $\mathcal{V}_{n+2,1}$, the details are left to the interested readers.    \vskip0.2cm

\subsection{Asymmetric valleys with weight $k+1$ in Dyck paths}

In this subsection, we follow with interest the left asymmetric valleys with weight $k+1$ in Dyck paths. The right asymmetric valleys are equivalent distribution to the left asymmetric valleys according to the symmetry of Dyck paths.

Let $\mathcal{V}^{L}_{n,k}$ denote the set of Dyck paths of length $2(n+3)$ having a distinguished left asymmetric valley with weight $k+1$. Set $V^{L}_{n,k}=|\mathcal{V}^{L}_{n,k}|$, which is the total number of left asymmetric valleys with weight $k+1$ in $\mathcal{D}_{n+3}$.

\begin{theorem}\label{theom 3.2.1}
There is a bijection between the sets $\mathcal{V}^{L}_{n,k}$ and $\mathcal{E}_{n+2,2k+2}$.
\end{theorem}

\pf Given a Dyck path $Q\in \mathcal{V}^{L}_{n,k}$ with a distinguished left asymmetric valley $\underline{\mathbf{d}^{k+j+2}\mathbf{u}^{k+1}}$ at level $i$ for certain $i, j\geq 0$, $Q$ can be uniquely partitioned into
\begin{eqnarray*}
Q=\left\{
\begin{array}{rl}
Q_0\mathbf{u}Q_1\dots \mathbf{u}Q_{k+2}\underline{\mathbf{d}^{k+2}\mathbf{u}^{k+1}}\mathbf{d}Q_{k+3}\dots\mathbf{d}Q_{2k+3}, & \mbox{when}\  i=0, \\[5pt]
Q_0\mathbf{u}Q_1\dots \mathbf{u}Q_{k+2}\underline{\mathbf{d}^{k+2}\mathbf{u}^{k+1}}\mathbf{d}Q_{k+3}\dots\mathbf{d}Q_{2k+3}\mathbf{d}{Q'}_0, & \mbox{when}\  i\geq 1,
\end{array}\right.
\end{eqnarray*}
where $Q_1, \dots, Q_{2k+3}$ are Dyck paths and $Q_{k+2}$ ends with $j$ $\mathbf{d}$ steps for certain $j\geq 0$ which, together with $\underline{\mathbf{d}^{k+2}\mathbf{u}^{k+1}}$, form the left asymmetric valley $\underline{\mathbf{d}^{k+j+2}\mathbf{u}^{k+1}}$ of $Q$, $Q_0$ is a partial Dyck path ending at level $i$, and $\overline{Q'}_0$ is a partial Dyck path ending at level $i-1\geq 0$.

In the $i=0$ case, $Q_0$ is always a Dyck path, we define $\rho(Q)=(Q_0, Q_1\mathbf{u}Q_{2}\dots \mathbf{u}Q_{2k+3})\in \mathcal{E}_{n+2,2k+2}$. In the $i\geq 1$ case, we define $\rho(Q)=(Q_0'\mathbf{d}Q_0, Q_1\mathbf{u}Q_{2}\dots \mathbf{u}Q_{2k+3})\in \mathcal{E}_{n+2,2k+2}$. Note that $Q_0$ always begins with a $\mathbf{u}$ step, and $Q_0'\mathbf{d}Q_0$ is a free Dyck path with a lowest valleys at the level $-i\leq -1$ such that the leftmost lowest valley is the intersection of $Q_0'\mathbf{d}$ and $Q_0$.

Similarly, one can verify that $\rho$ is a bijection between the sets $\mathcal{V}^{L}_{n,k}$ and $\mathcal{E}_{n+2,2k+2}$, the details are left to the interested readers.   \qed \vskip0.2cm

Let $V_k^{L}(x)=\sum_{n\geq k}V^{L}_{n,k}x^n$, then $V_k^{L}(x)=\frac{1}{x^2}E_{2k+3}(x)$ by Theorem \ref{theom 3.2.1}. Together with (\ref{eqn 2.3}), we have

\begin{corollary}\label{coro 3.2.2} The generating function
\begin{eqnarray*}
V_k^{L}(x)=\sum_{n\geq k}V^{L}_{n,k}x^n = \frac{x^{2k}C(x)^{2k+3}}{\sqrt{1-4x}},
\end{eqnarray*}
and the triangle $\mathbf{V}^{L}=\big(V^{L}_{n,k}\big)_{n\geq k\geq 0}$ forms a Riordan array $\big(\frac{C(x)^3}{\sqrt{1-4x}}, x^2C(x)^2\big)$ with the general entry
$$V^{L}_{n,k}=\binom{2n-2k+3}{n-2k}.$$
\end{corollary}
The first values of $V^{L}_{n,k}$ are displayed in Table 3.2.

\begin{center}
\begin{eqnarray*}
\begin{array}{|c|cccc|}\hline
n/k & 0    & 1      & 2    & 3              \\\hline
  0 & 1    &        &      &               \\
  1 & 5    &        &      &               \\
  2 & 21   & 1      &      &               \\
  3 & 84   & 7      &      &               \\
  4 & 330  & 36     & 1    &               \\
  5 & 1287 & 165    & 9    &                \\
  6 & 5005 & 715    & 55   &  1             \\ \hline
\end{array}
\end{eqnarray*}
Table 3.2. The first values of $V^{L}_{n,k}$.
\end{center}

Theorem \ref{theom 3.1.1} and \ref{theom 3.2.1} suggest the following result, whose direct bijective proof is left to the interested readers.

\begin{corollary}\label{coro 3.2.3}
There is a bijection between the sets $\mathcal{V}_{n+2,k+1}$ and $\mathcal{V}^{L}_{n,k}$.
\end{corollary}

Let $\mathcal{V}^{*}_{n,k}$ denote the set of Dyck paths of length $2(n+2)$ having a distinguished symmetric or left asymmetric valley with weight $k+1$. Set ${V}^{*}_{n,k}=|\mathcal{V}^{*}_{n,k}|$. Note that
$\mathcal{V}^{*}_{n,k}=\mathcal{V}_{n,k}\cup \mathcal{V}^{L}_{n-1,k}$ and ${V}^{*}_{n,k}=V_{n,k}+V^{L}_{n-1,k}$. By Theorem \ref{theom 3.1.1} and \ref{theom 3.2.1}, together with \ref{eqn 2.3}, we have
\begin{corollary}\label{coro 3.2.4}
The generating function
\begin{eqnarray*}
{V}^{*}_k(x)=\sum_{n\geq k}{V}^{*}_{n,k}x^n =V_k(x)+xV^{L}_k(x)=\frac{x^{2k}C(x)^{2k+2}}{\sqrt{1-4x}},
\end{eqnarray*}
and the triangle $\mathbf{V}^{*}=\big({V}^{*}_{n,k}\big)_{n\geq k\geq 0}$ forms a Riordan array $\big(\frac{C(x)^2}{\sqrt{1-4x}}, x^2C(x)^2\big)$ with the general entry
$${V}^{*}_{n,k}=\binom{2n-2k+2}{n-2k}.$$
\end{corollary}

The first values of ${V}^{*}_{n,k}$ are presented in Table 3.3.

\begin{center}
\begin{eqnarray*}
\begin{array}{|c|cccc|}\hline
n/k & 0    & 1    & 2    & 3             \\\hline
  0 & 1    &      &      &              \\
  1 & 4    &      &      &              \\
  2 & 15   & 1    &      &              \\
  3 & 56   & 6    &      &              \\
  4 & 210  & 28   & 1    &             \\
  5 & 729  & 120  & 8    &             \\
  6 & 3003 & 495  & 45   & 1            \\\hline
\end{array}
\end{eqnarray*}
Table 3.3. The first values of ${V}^{*}_{n,k}$.
\end{center}

\vskip0.5cm

\section{Symmetric and asymmetric peaks with weight $k+1$ in partial Dyck paths}

\subsection{Symmetric peaks with weight $k+1$ in partial Dyck paths}

In this subsection, we focus on the symmetric peaks with weight $k+1$ in partial Dyck paths.

Let $\mathcal{S}^{P}_{n,k,r}$ denote the set of partial Dyck paths in $\mathcal{D}_{n,r}$ having a distinguished symmetric peak with weight $k+1$. Set $S^{P}_{n,k,r}=|\mathcal{S}^{P}_{n,k,r}|$, which is the total number of symmetric peaks with weight $k+1$ in $\mathcal{D}_{n,r}$.


\begin{lemma}\label{lemma 4.1.2}
The total number $\alpha_{n,k}$ of symmetric peaks with weight $k+1$ in $\mathbf{u}\mathcal{D}_{n}$ is counted by the generating function $\sum_{n\geq 0}\alpha_{n,k}x^n=\frac{x^{k+2}C(x)}{\sqrt{1-4x}}$.
\end{lemma}
\pf Let $\mathbf{u}P\in \mathbf{u}\mathcal{D}_{n}$, where $P$ is a Dyck path in $\mathcal{D}_{n}$. Note that a symmetric peak $\underline{\mathbf{u}^{k+1}\mathbf{d}^{k+1}}$ of weight $k+1$ in $P$ is also a symmetric peak of weight $k+1$ in $\mathbf{u}P$ if it is not at the beginning of $P$. If $P$ starts with a symmetric peak $\underline{\mathbf{u}^{k+1}\mathbf{d}^{k+1}}$, that is $P=\underline{\mathbf{u}^{k+1}\mathbf{d}^{k+1}}P_1$, where $P_1\in \mathcal{D}_{n-k-1}$, the first symmetric peak of $P$ becomes an asymmetric peak $\underline{\mathbf{u}^{k+2}\mathbf{d}^{k+1}}$ in $\mathbf{u}\mathcal{D}_{n}$. Clearly, there are $C_{n-k-1}$ such kind of symmetric peaks in all $P\in \mathcal{D}_{n}$, which is counted by $x^{k+1}C(x)$. Hence, by Corollary \ref{coro 2.1.3}, the total number $\alpha_{n,k}$ of symmetric peaks with weight $k+1$ in $\mathbf{u}\mathcal{D}_{n}$ is counted by the generating function
$$\sum_{n\geq 0}\alpha_{n,k}x^n=xS_k(x)-x^{k+1}C(x)=\frac{x^{k+2}C(x)}{\sqrt{1-4x}}.  $$     \vskip0.2cm

\begin{theorem}\label{theom 4.1.3}
The total number $S^{P}_{n,k,r}$ of symmetric peaks with weight $k+1$ in $\mathcal{D}_{n,r}$ is counted by the generating function $x^{k+r+1}C(x)^{r+1}\Big(1+\frac{(r+1)x}{\sqrt{1-4x}}\Big)$. Namely,
$$S^{P}_{n,k,r}=\frac{(r+1)\big((n-k+1)(n-k-r)-2\big)}{\big(2(n-k)-r-2\big)\big(2(n-k)-r-1\big)}\binom{2(n-k)-r-1}{n-k}$$
for $n\geq k+r+1$.
\end{theorem}
\pf For any $P\in \bigcup_{n\geq 0}\mathcal{D}_{n,r}$, $P$ can be uniquely written as $P=P_0\mathbf{u}P_1\mathbf{u}P_2\dots \mathbf{u}P_r$, where $P_0, P_1, \dots, P_r$ are Dyck paths. By Corollary \ref{coro 2.1.3}, the total number of symmetric peaks with weight $k+1$ in all $P_0$'s is counted by $xS_k(x)$ and the total number of paths $P_1\mathbf{u}P_2\dots \mathbf{u}P_r\in \bigcup_{n\geq 0}\mathcal{D}_{n,r-1}$ is counted by $x^{r-1}C(x)^r$. By Lemma \ref{lemma 4.1.2}, the total number of symmetric peaks with weight $k+1$ in all $\mathbf{u}P_i$'s for $1\leq i\leq r$ is counted by $\frac{x^{k+2}C(x)}{\sqrt{1-4x}}$ and the total number of paths $P_0\mathbf{u}P_1\dots \mathbf{u}P_{i-1}\mathbf{u}P_{i+1}\dots \mathbf{u}P_r\in \bigcup_{n\geq 0}\mathcal{D}_{n,r-1}$ is counted by $x^{r-1}C(x)^r$. So the total number $S^{P}_{n,k,r}$ of symmetric peaks with weight $k+1$ in $\mathcal{D}_{n,r}$ is counted by
$$xS_k(x)x^{r}C(x)^r+rx^{r-1}C(x)^rx\frac{x^{k+2}C(x)}{\sqrt{1-4x}}=x^{k+r+1}C(x)^{r+1}\Big(1+\frac{(r+1)x}{\sqrt{1-4x}}\Big).$$
By (\ref{eqn 1.2}) and (\ref{eqn 2.1}), one can deduce that
\begin{eqnarray*}
S^{P}_{n,k,r} &=& [x^n]x^{k+r+1}C(x)^{r+1}\Big(1+\frac{(r+1)x}{\sqrt{1-4x}}\Big) \\
              &=& \frac{r+1}{n-k}\binom{2(n-k)-r-2}{n-k-1}+(r+1)\binom{2(n-k)-r-3}{n-k-1}   \\
              &=& \frac{(r+1)\big((n-k+1)(n-k-r)-2\big)}{\big(2(n-k)-r-2\big)\big(2(n-k)-r-1\big)}\binom{2(n-k)-r-1}{n-k}.
\end{eqnarray*}
This completes the proof.  \qed\vskip0.2cm

\subsection{Asymmetric peaks with weight $k+1$ in partial Dyck paths}

In this subsection, we take into account the left asymmetric peaks with weight $k+1$ in partial Dyck paths.

Let $\mathcal{L}^{P}_{n,k,r}$ denote the set of partial Dyck paths in $\mathcal{D}_{n,r}$ having a distinguished left asymmetric peak with weight $k+1$. Set $L^{P}_{n,k,r}=|\mathcal{L}^{P}_{n,k,r}|$, which is the total number of left asymmetric peaks with weight $k+1$ in $\mathcal{D}_{n,r}$.

\begin{lemma}\label{lemma 4.2.2}
The total number $\beta_{n,k}$ of left asymmetric peaks with weight $k+1$ in $\mathbf{u}\mathcal{D}_{n}$ is counted by the generating function $\sum_{n\geq 0}\beta_{n,k}x^n=\frac{x^{k+1}}{2}\Big(1+\frac{1}{\sqrt{1-4x}}\Big)=\frac{x^{k+1}}{C(x)\sqrt{1-4x}}$.
\end{lemma}
\pf Let $\mathbf{u}P\in \mathbf{u}\mathcal{D}_{n}$, where $P$ is a Dyck path in $\mathcal{D}_{n}$. Note that an asymmetric peak $\underline{\mathbf{u}^{k+j+2}\mathbf{d}^{k+1}}$ of weight $k+1$ in $P$ is also an asymmetric peak of weight $k+1$ in $\mathbf{u}P$. If $P$ starts with a symmetric peak $\underline{\mathbf{u}^{k+1}\mathbf{d}^{k+1}}$, that is $P=\underline{\mathbf{u}^{k+1}\mathbf{d}^{k+1}}P_1$, where $P_1\in \mathcal{D}_{n-k-1}$, the first symmetric peak of $P$ becomes an asymmetric peak $\underline{\mathbf{u}^{k+2}\mathbf{d}^{k+1}}$ in $\mathbf{u}P$. Clearly, there are $C_{n-k-1}$ such kind of symmetric peaks in all $P\in \mathcal{D}_{n}$, which is counted by $x^{k+1}C(x)$. Hence, by Corollary \ref{coro 2.2.3}, the total number $\beta_{n,k}$ of asymmetric peaks with weight $k+1$ in $\mathbf{u}\mathcal{D}_{n}$ is counted by
\begin{eqnarray*}
\sum_{n\geq 0}\beta_{n,k}x^n &=& x^3L_k(x)+x^{k+1}C(x) = x^{k+1}C(x)\Big(\frac{x^2C(x)^2}{\sqrt{1-4x}}+1\Big) \\
                             &=& \frac{x^{k+1}}{2}\Big(1+\frac{1}{\sqrt{1-4x}}\Big)=\frac{x^{k+1}}{C(x)\sqrt{1-4x}},
\end{eqnarray*}
where we use the relations $\sqrt{1-4x}=1-2xC(x)$ and $C(x)=1+xC(x)^2=\frac{1}{1-xC(x)}$. \qed\vskip0.2cm

\begin{theorem}\label{theom 4.2.3}
The total number $L^{P}_{n,k,r}$ of asymmetric peaks with weight $k+1$ in $\mathcal{D}_{n,r}$ is counted by the generating function $\frac{x^{k+r+1}C(x)^{r-1}}{\sqrt{1-4x}}\Big(r+x^2C(x)^4\Big)$. That is
$$L^{P}_{n,k,r}= r\binom{2(n-k)-r-3}{n-k-r-1}+\binom{2(n-k)-r-3}{n-k}$$
for $n\geq k+r+1$.
\end{theorem}
\pf For any $P\in \bigcup_{n\geq 0}\mathcal{D}_{n,r}$, $P$ can be uniquely written as $P=P_0\mathbf{u}P_1\mathbf{u}P_2\dots \mathbf{u}P_r$, where $P_0, P_1, \dots, P_r$ are any Dyck paths. By Corollary \ref{coro 2.2.3}, the total number of asymmetric peaks with weight $k+1$ in all $P_0$'s is counted by $x^3L_k(x)$ and the total number of paths $P_1\mathbf{u}P_2\dots \mathbf{u}P_r\in \bigcup_{n\geq 0}\mathcal{D}_{n,r-1}$ is counted by $x^{r-1}C(x)^r$. By Lemma \ref{lemma 4.2.2}, the total number of asymmetric peaks with weight $k+1$ in all $\mathbf{u}P_i$'s for $1\leq i\leq r$ is counted by $\frac{x^{k+1}}{C(x)\sqrt{1-4x}}$ and the total number of paths $P_0\mathbf{u}P_1\dots \mathbf{u}P_{i-1}\mathbf{u}P_{i+1}\dots \mathbf{u}P_r\in \bigcup_{n\geq 0}\mathcal{D}_{n,r-1}$ is counted by $x^{r-1}C(x)^r$. So the total number $L^{P}_{n,k,r}$ of symmetric peaks with weight $k+1$ in $\mathcal{D}_{n,r}$ is counted by
$$x^3L_k(x)x^{r}C(x)^r+rx^{r-1}C(x)^rx\frac{x^{k+1}}{C(x)\sqrt{1-4x}}=\frac{x^{k+r+1}C(x)^{r-1}}{\sqrt{1-4x}}\Big(r+x^2C(x)^4\Big).$$
By (\ref{eqn 1.2}) and (\ref{eqn 2.1}), one can deduce that
\begin{eqnarray*}
L^{P}_{n,k,r} &=& [x^n]\frac{x^{k+r+1}C(x)^{r-1}}{\sqrt{1-4x}}\Big(r+x^2C(x)^4\Big) \\
              &=& r\binom{2(n-k)-r-3}{n-k-r-1}+\binom{2(n-k)-r-3}{n-k}.
\end{eqnarray*}
This completes the proof.  \qed

\begin{remark}
One can also discuss the total number of symmetric peaks or left asymmetric peaks in all primitive Dyck paths of a given length. By the similar methods, one can deduce that the total number of symmetric peaks of weight $k+1$
in all primitive Dyck paths of length $2(n+4)$ is counted by the generating function $\frac{x^{k}C(x)^3}{\sqrt{1-4x}}$, and the total number of left asymmetric peaks of weight $k+1$
in all primitive Dyck paths of length $2(n+3)$ is counted by the generating function $\frac{x^{k}C(x)}{\sqrt{1-4x}}$. Hence, there exist bijections between the set of Dyck paths of length $2(n+3)$ having a distinguished left asymmetric peak of weight $k+1$ and the set of primitive Dyck paths of length $2(n+4)$ having a distinguished symmetric peak of weight $k+1$, and bijections between the set of Dyck paths of length $2(n+2)$ having a distinguished symmetric valley of weight $1$ and the set of primitive Dyck paths of length $2(n+3)$ having a distinguished left asymmetric peak of weight $1$. The details are left to interested readers.
\end{remark}

\vskip0.5cm
\section*{Declaration of competing interest}

The authors declare that they have no known competing financial interests or personal relationships that could have
appeared to influence the work reported in this paper.

\section*{Acknowledgements} {The authors are grateful to the referees for
the helpful suggestions and comments. The Project is sponsored by ``Liaoning
BaiQianWan Talents Program". }

\vskip.2cm


\end{document}